\pdfoutput=1
\RequirePackage{ifpdf}
\ifpdf % We are running pdfTeX in pdf mode
\documentclass[pdftex]{sigma}
\else
\documentclass{sigma}
\fi

\begin{document}

%\allowdisplaybreaks

\renewcommand{\PaperNumber}{059}

\FirstPageHeading

\ShortArticleName{Formal Integrability for the Inverse Problem of the Calculus of Variations}

\ArticleName{Formal Integrability for the Nonautonomous Case\\ of the Inverse Problem of the Calculus of Variations}

\Author{Oana CONSTANTINESCU}

\AuthorNameForHeading{O.~Constantinescu}

\Address{Faculty of Mathematics, Alexandru Ioan Cuza University,\\
Bd.~Carol no.~11, 700506, Iasi, Romania}
\Email{\href{mailto:oanacon@uaic.ro}{oanacon@uaic.ro}}
\URLaddress{\url{http://www.math.uaic.ro/~oanacon/}}

\ArticleDates{Received March 16, 2012, in f\/inal form September 03, 2012; Published online September 06, 2012}

\Abstract{We address the integrability conditions of the inverse problem of
the calculus of variations for time-dependent SODE using the Spencer
version of the Cartan--K\"ahler theorem. We consider a linear partial
dif\/ferential operator $P$ given by the two Helmholtz conditions expressed
in terms of semi-basic 1-forms and study its formal integrability.
We prove that~$P$ is involutive and there is only one obstruction
for the formal integrability of this operator. The obstruction is
expressed in terms of the curvature tensor~$R$ of the induced nonlinear
connection. We recover some of the classes of Lagrangian semisprays:
f\/lat semisprays, isotropic semisprays and arbitrary semisprays on
2-dimensional manifolds.}

\Classification{49N45; 58E30; 34A26; 37J30}

\section{Introduction}

One of the most interesting problems of geometric mechanics is related
to the integrability conditions of the inverse problem of the calculus
of variations for time-dependent second-order ordinary dif\/ferential
equations (SODE). The inverse problem can be formulated as follows. Given
a~time-dependent system of SODE
\[
\frac{d^{2}x^{i}}{dt^{2}}+2G^{i}\left(t,x,\frac{dx}{dt}\right)=0,\qquad i\in \{1,\dots,n\},
\]
 under what conditions this system can be made equivalent, using a
multiplier matrix $g_{ij}$, with the system of Euler--Lagrange equations
of a regular Lagrangian \[
g_{ij}\left(t,x,\frac{dx}{dt}\right)\left(\frac{d^{2}x^{i}}{dt^{2}}+2G^{i}\left(t,x,\frac{dx}{dt}\right)\right)=\frac{d}{dt}\left(\frac{\partial L}{\partial y^{i}}\right)-\frac{\partial L}{\partial x^{i}}\,\,?
\]

In this case such a system is called variational. The necessary and
suf\/f\/icient conditions under which such a system is variational are
known as the \emph{Helmholtz conditions}.

This inverse problem was solved for the case $n=1$ by Darboux~\cite{Dar},
and for $n=2$ by Douglas~\cite{douglas41}. Douglas's approach consists
in an application of the Riquier theory of systems of partial dif\/ferential
equations~\cite{Rhi}, to a certain associated linear dif\/ferential
system. The generalization of its results in the higher dimensional
case is a very dif\/f\/icult problem because the system provided by the
Helmholtz conditions is extremely over-determined. Some of the f\/irst
studies of the inverse problem in spaces of arbitrary dimension are
those of Davis~\cite{davis29} and Kosambi~\cite{kosambi35}.

%\looseness=-1
There are dif\/ferent attempts to solve this problem. First, there are
some reformulations of the Helmholtz conditions in better geometric
forms, which are close enough to the f\/irst analytical formulations
\cite{douglas41,santilli78,Sarlet81,sarlet82}, but undercover more
of the geometry behind them \cite{carinena89,crampin81,crampin84a,crampin94,deleon88,krupkova08,massa94}.
The system of SODE is identif\/ied with a semispray on the f\/irst jet
bundle of a f\/ibred manifold over~$\mathbb{R}$. The most important
geometric tools induced by a semispray are nonlinear connection, Jacobi
endomorphism, dynamical covariant derivative, linear connections and
their curvatures. Some reformulations of the Helmholtz conditions
are using either the special derivations along the tangent bundle
projection introduced in \cite{sarlet95}, or the semi-basic 1-forms~\cite{B-C-2010} and the Fr\"olicher--Nijenhuis theory of derivations
on the algebra of vector-valued forms~\cite{frolicher56}.

Anderson and Thompson \cite{anderso92} analyzed the inverse problem
based on the exterior dif\/ferential system approach \cite{BCG91}.
Using the variational bicomplex associated to a system of arbitrary
order ordinary dif\/ferential equations, they derived the fundamental
system of equations for the va\-riatio\-nal multiplier and proved their
suf\/f\/iciency. They made a detailed study of two dimensional sprays
and they proved, for general degrees of freedom, that all isotropic
semisprays are va\-ria\-tional. It means semisprays that have the associated
Jacobi endomorphism a multiple of the identity. This correspond to
the \emph{Case~I} of Douglas's classif\/ication. This approach is continued
in~\cite{Ald-Pr-Sa-TH}, where the case of $\Phi$ diagonalizable,
with distinct eigenfunctions, is exposed in detail. The same Case
I was proved to be variational also in~\cite{Sar-Cram-Mart-integr}.
This paper uses Riquier theory, but in a more geometric way. The process
of repeated dif\/ferentiations of equations and searching for new nontrivial
relations is realized by intrinsic operations.

\looseness=1
Another subcase of Douglas's case II is discussed in \cite{Cr-Pr-Sa-Th-sep}:
separable systems of SODE. Any systems of SODE from this subcase is
variational. They showed that any system of SODE in~\emph{Case II1}
with $n$ degrees of freedom can be separated into n separate systems
of two f\/irst-order equations. They also proved that there are systems
separable in the above sense but not separable into single independent
second-order equations. This case was treated in~\cite{cantrijn96}.
In~\cite{sarlet02} the authors reinvestigated the case $n=2$ with
their more intrinsic version of the Riquier algorithm. Their approach
is based on the same underlying methodology as the analytical work
of Douglas.

Another method of studying the integrability conditions of the inverse
problem of the calculus of variation is the Spencer--Goldsmchmidt theory
of formal integrability of partial dif\/ferential operators, using two
suf\/f\/icient conditions provided by Cartan--K\"ahler theorem~\cite{Cartan,Gold,Spencer}.
This method was applied for autonomous SODE in~\cite{grifone00},
using the Fr\"olicher--Nijenhuis theory of derivations of vector-valued
dif\/ferential forms. Grifone and Muzsnay gave the f\/irst obstructions
so that a \emph{spray} (homogeneous semispray) is variational, for
general degrees of freedom. In order to obtain a complete classif\/ication
of variational sprays, they restricted their work to some particular
cases. The Spencer theory is fully applied to the two dimensional
case, corresponding to Douglas's paper. For the general n-dimensional
case, it is proved only that isotropic sprays are variational. It
is important to notice that Grifone and Muzsnay's analysis starts
from the Euler--Lagrange partial dif\/ferential operator, and not from
the Helmholtz conditions.

For time independent, homogeneous SODE, the inverse problem is known
as the projective metrizability problem. This problem and its formal
integrability is studied in~\cite{B-M_arxiv} using Spencer theory.
It was shown that there exists only one f\/irst obstruction for the
formal integrability of the projective metrizability operator, expressed
in terms of the curvature tensor of the nonlinear connection induced
by the spray. This obstruction correspond to second obstruction for
the formal integrability of the Euler--Lagrange operator.

An interesting and new approach regarding variational PDE's is the
one of A.~Pr\'asta\-ro~\mbox{\cite{Prast.99,Prast.06}}. Using suitable cohomologies
and integral bordism groups, the author characterizes variational
systems constrained by means of PDE's of submanifolds of f\/iber bundles.
He presents a~new algebraic topological characterization of global
solutions of variational problems.

In this paper we address the integrability conditions of the inverse
problem of the calculus of variations for time-dependent SODE using
also the Spencer version of the Cartan--K\"ahler theorem. The proper
setting is the f\/irst jet bundle $J^{1}\pi$ of an $(n+1)$ manifold~$M$ f\/ibred over~$\mathbb{R}$. In \cite{B-C-2010} it is proved that
a time-dependent semispray is Lagrangian if and only if there exists
a~semi-basic 1-form $\theta$ on~$J^{1}\pi$, that satisf\/ies a dif\/ferential
system. This gives rise to a linear partial dif\/ferential operator~$P$.
We study the formal integrability of~$P$ using two suf\/f\/icient conditions
provided by Cartan--K\"ahler theorem. We prove that the symbol $\sigma^{1}(P)$
is involutive (Theorem~\ref{thm:-involutivity}) and hence there is
only one obstruction for the formal integrability of the ope\-ra\-tor~$P$, which is due to curvature tensor~$R$ (Theorem~\ref{thm:A-first-order}).
Based on this result, we recover some of the classes of Lagrangian
semisprays: f\/lat semisprays, isotropic semisprays and arbitrary semisprays
on $2$-dimensional jet spaces ($n=1$).

The motivation for this article is double-folded. So far all the results
about the inverse problem of the calculus of variations were obtained
separately, in the autonomous and non\-auto\-nomous settings. This is
due to the dif\/ferent frameworks involved: the tangent bundle~$TM$
(a~vector bundle) and respectively the f\/irst jet bundle~$J^{1}\pi$
(an af\/f\/ine bundle). The geometric tools are usually constructed in
dif\/ferent ways, and special attention was given to the time-depending
situation. This paper follows the line of~\cite{B-M_arxiv} but naturally
the proofs of the main theorems have some particularities due to the
dif\/ferent setting.

Secondly, there are similarities between the formulation of the Helmholtz
conditions for sprays in the autonomous setting and respectively for
semisprays in the nonautonomous one \cite{B-C-2010, bucataru09}. This
is natural because $J^{1}\pi$ can be embedded in $\widetilde{TM}$
(the tangent bundle with the zero section removed). Due to this embedding
one can associate to any regular Lagrangian on the velocity-phase
space $J^{1}\pi$ a homogeneous degenerate Lagrangian on the extended
phase space $\widetilde{TM}$, such that the action def\/ined by a curve
in the jet formalism coincides with the action def\/ined by the corresponding
curve in the extended formalism. There are correspondences between
the main geometric objects associated to these Lagrangians: Poincar\'e
1-~and~2-forms, energies, canonical semisprays-sprays
\cite{Ant-B,CantrijnHom86,CarinenaHom2006,Klein}. Therefore,
due to this homogeneous formalism, it is natural to expect such kind
of similarities between the results corresponding to homogeneous structures
on~$\widetilde{TM}$ and nonhomogeneous one on~$J^{1}\pi$.

The paper is organized as follows. In Section~\ref{sec:prelim} we introduce the principal
geometric tools induced by a time-dependent semispray on $J^{1}\pi$
and characterize Lagrangian vector f\/ields with respect to semi-basic
1-forms. Section~\ref{section3} is dedicated to the application of the
Spencer theory to the study of formal integrability of the partial dif\/ferential operator (PDO) $P=(d_{J},d_{h})$.
The most important results are Theorems \ref{thm:-involutivity} and~\ref{thm:A-first-order}. Section~\ref{section3.3} presents classes
of semisprays for which the obstruction in Theorem~\ref{thm:A-first-order}
is automatically satisf\/ied. For these classes, the PDO~$P$ is formally
integrable, and hence these semisprays will be Lagrangian SODE.

\section{Preliminaries} \label{sec:prelim}

\subsection[The first-order jet bundle $J^{1}\pi$]{The f\/irst-order jet bundle $\boldsymbol{J^{1}\pi}$} \label{subsec:fn}

The appropriate geometric setting for the study of time-dependent
SODE is the af\/f\/ine jet bundle $(J^{1}\pi,\pi_{10},M)$~\cite{saunders89}.
We consider an $(n+1)$-dimensional, real, smooth manifold $M$, which
is f\/ibred over $\mathbb{R}$, $\pi:M\rightarrow\mathbb{R}$, and represents
the space-time. The f\/irst jet bundle of $\pi$ is denoted by $\pi_{10}:J^{1}\pi\to M$,
$\pi_{10}(j_{t}^{1}\phi)=\phi(t),$ for $\phi$ a local section of
$\pi$ and $j_{t}^{1}\phi$ the f\/irst-order jet of~$\phi$ at~$t$.
A local coordinate system $(t,x^{i})_{i\in \{1,\dots,n\}}$ on $M$
induces a local coordinate system on $J^{1}\pi$, denoted by $(t,x^{i},y^{i})$.
Submersion $\pi_{10}$ induces a \emph{natural foliation} on $J^{1}\pi$
such that~$(t,x^{i})$ are transverse coordinates for this foliation,
while~$(y^{i})$ are coordinates for the leaves of the foliation.
Throughout the paper we consider Latin indices $i\in \{1,\dots,n\}$
and Greek indices $\alpha\in \{0,\dots,n\}$, using the notation~$\left(x^{\alpha}\right)=(t=x^{0},x^{i})$.

In this article we use the Fr\"olicher--Nijenhuis theory \cite{frolicher56,grifone00,KMS93}
of derivations of vector-valued dif\/ferential forms on the f\/irst jet
bundle $J^{1}\pi$. We adopt the following notations: $C^{\infty}(J^{1}\pi)$
for the ring of smooth functions on~$J^{1}\pi$, $\mathfrak{X}(J^{1}\pi)$
for the~$C^{\infty}$ module of vector f\/ields on $J^{1}\pi$ and $\Lambda^{k}(J^{1}\pi)$
for the~$C^{\infty}$ module of $k$-forms on~$J^{1}\pi$. The $C^{\infty}$
module of $(r,s)$-type vector f\/ields on~$J^{1}\pi$ is denoted by
$\mathcal{T}_{s}^{r}(J^{1}\pi)$ and the tensor algebra on~$J^{1}\pi$
is denoted by $\mathcal{T}(J^{1}\pi)$. The graded algebra of dif\/ferential
forms on $J^{1}\pi$ is written as $\Lambda(J^{1}\pi)=\bigoplus_{k\in \{1,\dots,2n+1\}}\Lambda^{k}(J^{1}\pi)$.
We denote by~$S^{k}(J^{1}\pi)$ the space of symmetric $(0,k)$ tensors
on~$J^{1}\pi$ and by $\Psi(J^{1}\pi)=\bigoplus_{k\in \{1,\dots,2n+1\}}\Psi^{k}(J^{1}\pi)$
the graded algebra of vector-valued dif\/ferential forms on $J^{1}\pi$.
Throughout the paper we assume that all objects are $C^{\infty}$-smooth
where def\/ined.

A parametrized curve on $M$ is a section of $\pi$: $\gamma:\mathbb{R}\rightarrow M$,
$\gamma(t)=(t,x^{i}(t))$. Its \emph{first-order jet prolongation}
$J^{1}\gamma:t\in\mathbb{R}\rightarrow J^{1}\gamma(t)=\left(t,x^{i}(t),dx^{i}/dt\right)\in J^{1}\pi$
is a section of the f\/ibration $\pi_{1}:=\pi\circ\pi_{10}:  J^{1}\pi\rightarrow\mathbb{R}$.

Let $VJ^{1}\pi$ be the \emph{vertical subbundle} of $TJ^{1}\pi,$
$VJ^{1}\pi=\{\xi\in TJ^{1}\pi,\,  D\pi_{10}(\xi)=0\}\subset TJ^{1}\pi$.
The f\/ibers $V_{u}J^{1}\pi=\operatorname{Ker} D_{u}\pi_{10}$, $u\in J^{1}\pi$ determine
a regular, $n$-dimensional, integrable vertical distribution. Remark
that $V_{u}J^{1}\pi=\operatorname{spann}\{\partial/\partial y^{i}\}$
and its annihilators are the \emph{contact $1$-forms} $\delta x^{i}=dx^{i}-y^{i}dt$, $i\in \{1,\dots,n\}$
and \emph{basic $1$-forms} $\lambda dt$, $\lambda\in C^{\infty}(J^{1}\pi)$.
The \emph{vertical endomorphism} $J=\frac{\partial}{\partial y^{i}}\otimes\delta x^{i}$
is a vector-valued 1-form on~$J^{1}\pi$, with $\operatorname{Im}J=V(J^{1}\pi)$,
$V(J^{1}\pi)\subset\operatorname{Ker}J$ and $J^{2}=0$.

Its Fr\"olicher--Nijenhuis tensor is given by
\begin{gather*}
N_{J}=\frac{1}{2}[J,J]=-\frac{\partial}{\partial y^{i}}\otimes\delta x^{i}\wedge dt=-J\wedge dt.
%\label{eq:NijenhuisJ}
\end{gather*}
 Consequently, $d_{J}^{2}=d_{N_{J}}=-d_{J\wedge dt}\neq0$ and therefore
$d_{J}$-exact forms on $J^{1}\pi$ may not be $d_{J}$-closed. Here
$d_{J}$ is the \emph{exterior derivative} with respect to the vertical
endomorphisms.

\begin{remark}
For $A\in\Psi^{1}(J^{1}\pi)$ a vector-valued 1-form, the exterior
derivative with respect to $A$ is a derivation of degree 1 given
by $d_{A}=i_{A}\circ d-d\circ i_{A}.$
\end{remark}
A $k$-form $\omega$ on $J^{1}\pi$, $k\geq1$, is called \emph{semi-basic}
if it vanishes whenever one of the arguments is vertical.

A vector-valued $k$-form $A$ on $J^{1}\pi$ is called semi-basic
if it takes values in the vertical bundle and it vanishes whenever
one of the arguments is vertical.

A semi-basic $k$-form satisf\/ies the relation $i_{J}\theta=0$ and
locally can be expressed as $\theta=\theta_{0}dt+\theta_{i}\delta x^{i}$.
For example, contact 1-forms $\delta x^{i}$ are semi-basic 1-forms.

If a vector-valued $k$-form $A$ is semi-basic, then $J\circ A=0$
and $i_{J}A=0$. The vertical endomorphism $J$ is a vector-valued,
semi-basic $1$ -form.

Locally, a semi-basic $k$-form $\theta$ has the next form \begin{eqnarray*}
\theta=\frac{1}{k!}\theta_{i_{1}\dots i_{k}}(x^{\alpha},y^{j})\delta x^{i_{1}}\wedge\cdots\wedge\delta x^{i_{k}}+\frac{1}{(k-1)!}\widetilde{\theta}_{i_{1}\dots i_{k-1}}(x^{\alpha},y^{j})\delta x^{i_{1}}\wedge\cdots\wedge\delta x^{i_{k-1}}\wedge dt.\end{eqnarray*}

For simplicity, we denote by $T^{*}$ the vector bundle of $1$-forms
on $J^{1}\pi$, by $T_{v}^{*}$ the vector bundle of semi-basic $1$-forms
on $J^{1}\pi$ and by $\Lambda^{k}T_{v}^{*}$ the vector bundle of
semi-basic $k$-forms on $J^{1}\pi$. We also denote by $\Lambda_{v}^{k}=\operatorname{Sec}\left(\Lambda^{k}T_{v}^{*}\right)$
the $C^{\infty}(J^{1}\pi)$-module of sections of $\Lambda^{k}T_{v}^{*}$
and by $S^{k}T^{*}$ the vector bundle of symmetric tensors of $(0,k)$-type
on $J^{1}\pi$. $S^{1}T^{*}$ will be identif\/ied with $T^{*}$.

A \emph{semispray} is a globally def\/ined vector f\/ield $S$ on $J^{1}\pi$
such that \begin{gather*}
J(S)=0\qquad\textrm{and}\qquad dt(S)=1.%\label{eq:4}
\end{gather*}
 The integral curves of a semispray are f\/irst-order jet prolongations
of sections of $\pi\circ\pi_{10}:\!  J^{1}\pi\!\rightarrow\!\mathbb{R}$.
Locally, a semispray has the form{\samepage \begin{gather}
S=\frac{\partial}{\partial t}+y^{i}\frac{\partial}{\partial x^{i}}-2G^{i}(x^{\alpha},y^{j})\frac{\partial}{\partial y^{i}},\label{eq:5}
\end{gather}
 where functions $G^{i}$, called the semispray coef\/f\/icients, are
locally def\/ined on~$J^{1}\pi$.}

A parametrized curve $\gamma:  I\rightarrow M$ is a \emph{geodesic}
of $S$ if $S\circ J^{1}\gamma=\frac{d}{dt}(J^{1}\gamma).$

In local coordinates, $\gamma(t)=(t,x^{i}(t))$ is a geodesic of the
semispray $S$ given by \eqref{eq:5} if and only if it satisf\/ies
the system of SODE \begin{gather}
\frac{d^{2}x^{i}}{dt^{2}}+2G^{i}\left(t,x,\frac{dx}{dt}\right)=0.\label{sode}
\end{gather}

Therefore such a system of time-dependent SODE can be identif\/ied with
a semispray on $J^{1}\pi$.

{\bf Canonical nonlinear connection.}
A~\emph{nonlinear connection} on $J^{1}\pi$ is an $(n+1)$-dimensional
distribution $H:u\in J^{1}\pi\mapsto H_{u}\subset T_{u}J^{1}\pi$,
supplementary to $VJ^{1}\pi$: $\forall\, u\in J^{1}\pi$, $T_{u}J^{1}\pi=H_{u}\oplus V_{u}$.

A semispray $S$ induces a nonlinear connection on $J^{1}\pi$, given
by the \emph{almost product structure} $\Gamma=-\mathcal{L}_{S}J+S\otimes dt$, $\Gamma^{2}={\rm Id}$.
The \emph{horizontal projector} that corresponds to this almost product
structure is $h=\frac{1}{2}\left({\rm Id}-\mathcal{L}_{S}J+S\otimes dt\right)$
and the \emph{vertical projector} is $v={\rm Id}-h.$

The horizontal subspace is spanned by $S$ and by $\frac{\delta}{\delta x^{i}}:=\frac{\partial}{\partial x^{i}}-N_{i}^{j}\frac{\partial}{\partial y^{j}}$, where $N_{j}^{i}=\frac{\partial G^{i}}{\partial y^{j}}$.
In this paper we prefer to work with the following adapted basis and
cobasis:
\begin{gather}
\left\{ S,\frac{\delta}{\delta x^{i}},\frac{\partial}{\partial y^{i}}\right\} ,\qquad\{dt,\delta x^{i},\delta y^{i}\},\label{eq:11}
\end{gather}
 with $\delta x^{i}$ the \emph{contact $1$-forms} and $\delta y^{i}=dy^{i}+N_{\alpha}^{i}dx^{\alpha}$, $N_{0}^{i}=2G^{i}-N_{j}^{i}y^{j}$.
Functions $N_{j}^{i}$ and $N_{0}^{i}$ are the coef\/f\/icients of the
nonlinear connection induced by the semispray~$S$.

With respect to basis and cobasis \eqref{eq:11}, the horizontal and
vertical projectors are locally expressed as $h=S\otimes dt+\frac{\delta}{\delta x^{i}}\otimes\delta x^{i}$, $v=\frac{\partial}{\partial y^{i}}\otimes\delta y^{i}$.
We consider the $(1,1)$-type tensor f\/ield $\mathbb{F}=h\circ\mathcal{L}_{S}h-J$,
which corresponds to the almost complex structure in the autonomous
case. It satisf\/ies $\mathbb{F}^{3}+\mathbb{F}=0,$ which means that
it is an~$f(3,1)$ structure. It can be expressed locally as $\mathbb{F}=\frac{\delta}{\delta x^{i}}\otimes\delta y^{i}-\frac{\partial}{\partial y^{j}}\otimes\delta x^{i}$.

{\bf Curvature.}
The following properties for the torsion and curvature of the nonlinear
connection induced by the semispray are proved in~\cite{B-C-2010}.

The weak torsion tensor f\/ield of the nonlinear connection $\Gamma$
vanishes: $[J,h]=0,$ which is equivalent also with $[J,\Gamma]=0$.

The curvature tensor $R=N_{h}$ of the nonlinear connection $\Gamma$
is a vector-valued semi-basic 2-form, locally given by
\begin{gather}
R=\frac{1}{2}[h,h]=\frac{1}{2}R_{ij}^{k}\frac{\partial}{\partial y^{k}}\otimes\delta x^{i}\wedge\delta x^{j}+R_{i}^{j}\frac{\partial}{\partial y^{j}}\otimes dt\wedge\delta x^{i},\label{eq:19}
\end{gather}
 where \begin{gather*}
R_{jk}^{i}=\frac{\delta N_{j}^{i}}{\delta x^{k}}-\frac{\delta N_{k}^{i}}{\delta x^{j}}%\label{eq:Rijk}
\end{gather*}
 and \begin{gather}
R_{j}^{i}=2\frac{\partial G^{i}}{\partial x^{j}}-\frac{\partial G^{i}}{\partial y^{k}}\frac{\partial G^{k}}{\partial y^{j}}-S\left(\frac{\partial G^{i}}{\partial y^{j}}\right).\label{eq:local-Jacobi}
\end{gather}

The \emph{Jacobi endomorphism} is def\/ined as \begin{gather}
\Phi=v\circ\mathcal{L}_{S}h=\mathcal{L}_{S}h-\mathbb{F}-J.\label{eq:22}\end{gather}
 Jacobi endomorphism $\Phi$ is a semi-basic, vector-valued 1-form
and satisf\/ies $\Phi^{2}=0$. Locally, can be expressed as $\Phi=R_{i}^{j}\frac{\partial}{\partial y^{j}}\otimes\delta x^{i},$
where $R_{j}^{i}$ are given by~(\ref{eq:local-Jacobi}).

The Jacobi endomorphism and the curvature of the nonlinear connection
are related by the following formulae:
\begin{gather}
\Phi   =   i_{S}R,\label{eq:23'}\\
\left[J,\Phi\right]   =   3R+\Phi\wedge dt.\label{eq:24'}
\end{gather}

Remark that $R=0$ if and only if $\Phi=0$.

\begin{definition}
A semispray $S$ is called \emph{isotropic} if its Jacobi endomorphism
has the form \begin{gather}
\Phi=\lambda J,
\label{eq:phi_iso}
\end{gather}
 where $\lambda\in C^{\infty}(J^{1}\pi)$.
\end{definition}
Next we express the isotropy condition \eqref{eq:phi_iso} for a semispray
in terms of the curvature tensor~$R$.

\begin{proposition}
\label{prop:iso} A semispray $S$ is isotropic if and only if its
curvature tensor $R$ has the form
\begin{gather*}
R=\alpha\wedge J,%\label{eq:r_iso}
\end{gather*}
 where $\alpha$ is a semi-basic $1$-form on $J^{1}\pi$.
\end{proposition}

\begin{proof}
Suppose that $S$ is an isotropic SODE. Then there exists $\lambda\in C^{\infty}(J^{1}\pi)$
such that $\Phi=\lambda J$. From (\ref{eq:24'}) it results
\begin{gather*}
3R   =   [J,\lambda J]-\Phi\wedge dt,\\
{}[J,\lambda J]   =   \left(d_{J}\lambda\right)\wedge J-d\lambda\wedge J^{2}+\lambda[J,J] \ \Rightarrow \
R   =   \frac{1}{3}\left(d_{J}\lambda\right)\wedge J-\lambda J\wedge dt=\alpha\wedge J,
\end{gather*}
 with $\alpha=\frac{1}{3}d_{J}\lambda+\lambda dt\in T_{v}^{*}$.

In the above calculus we used the formula \cite{grifone00}\[
[K,gL]= (d_{K}g )\wedge L-dg\wedge KL+g[K,L],\]
 for $K$, $L$ vector-valued one-forms on $J^{1}\pi$ and $g\in C^{\infty}(J^{1}\pi)$.

For the converse, suppose that $R=\alpha\wedge J$, with $\alpha\in T_{v}^{*}$.
Formula (\ref{eq:23'}) implies $\Phi=i_{S} (\alpha\wedge J )= (i_{S}\alpha )J-\alpha\wedge i_{S}J= (i_{S}\alpha )J$.
\end{proof}

\subsection{Lagrangian semisprays}

In this subsection we recall some basic notions about Lagrangian semisprays.

\begin{definition}\null \qquad

1) A smooth function $L\in C^{\infty}(J^{1}\pi)$ is called a \emph{Lagrangian
function}.

2) The Lagrangian $L$ is regular if the $(0,2)$ type tensor with
local components
\begin{gather*}
g_{ij}\big(x^{\alpha},y^{k}\big)=\frac{\partial^{2}L}{\partial y^{i}\partial y^{j}}%\label{eq:tensorg_ij}
\end{gather*}
has rank $n$ on $J^{1}\pi$. The tensor $g=g_{ij}\delta x^{i}\otimes\delta x^{j}$
is called the \emph{metric tensor} of the Lagrangian~$L$.
\end{definition}

\begin{remark}
More exactly, \cite{saunders89}, a function $L\in C^{\infty}(J^{1}\pi)$
is called a Lagrangian density on~$\pi$. If~$\Omega$ is a volume
form on $\mathbb{R}$, the corresponding Lagrangian is the semi-basic
1-form $L\pi_{1}^{*}\Omega$ on~$J^{1}\pi$. Using a f\/ixed volume form
on~$\mathbb{R}$, for example~$dt$, it is natural to consider the
function~$L$ as a~(f\/irst-order) Lagrangian.
\end{remark}

For the particular choice of $dt$ as volume form on $\mathbb{R}$,
the \emph{Poincar\'e--Cartan $1$-form} of the Lagrangian $L$ is
$\theta_{L}:=Ldt+d_{J}L$. The Lagrangian $L$ is regular if and only
if the \emph{Poincar\'e--Cartan $2$-form} $d\theta_{L}$ has maximal
rank $2n$ on $J^{1}\pi$.

For a detailed exposition on the regularity conditions for Lagrangians
see~\cite{Krupkova}.

The \emph{geodesics} of a semispray $S$, given by the system of SODE
(\ref{sode}), coincide with the solutions of the \emph{Euler--Lagrange
equations} \begin{gather*}
\frac{d}{dt}\left(\frac{\partial L}{\partial y^{i}}\right)-\frac{\partial L}{\partial x^{i}}=0%\label{eq:EL}
\end{gather*}
 if and only if \begin{gather}
g_{ij}\left(t,x,\frac{dx}{dt}\right)\left(\frac{d^{2}x^{i}}{dt^{2}}+2G^{i}\left(t,x,\frac{dx}{dt}\right)\right)=\frac{d}{dt}\left(\frac{\partial L}{\partial y^{i}}\right)-\frac{\partial L}{\partial x^{i}}.\label{eq:SEL}\end{gather}

Therefore, for a semispray $S,$ there exists a Lagrangian function
$L$ such that (\ref{eq:SEL}) holds true if and only if $S\left(\frac{\partial L}{\partial y^{i}}\right)-\frac{\partial L}{\partial x^{i}}=0,$
which can be further expressed as
\begin{gather}
\mathcal{L}_{S}\theta_{L}=dL \ \Leftrightarrow \ i_{S}d\theta_{L}=0.\label{eq:LSTL}
\end{gather}

\begin{definition}
A semispray $S$ is called a \emph{Lagrangian semispray} (or a Lagrangian
vector f\/ield) if there exists a Lagrangian function~$L$,
locally def\/ined on~$J^{1}\pi$, that satisf\/ies~(\ref{eq:LSTL}).
\end{definition}

In \cite{B-C-2010} it has been shown that a semispray $S$ is a Lagrangian
semispray if and only if there exists a semi-basic 1-form $\theta\in\Lambda_{v}^{1}$
with $\operatorname{rank}(d\theta)=2n$ on $J^{1}\pi$, such that
$\mathcal{L}_{S}\theta$ is closed. This represents a reformulation,
in terms of semi basic $1$-forms, of the result in terms of $2$-forms
obtained by Crampin et al.\ in~\cite{crampin84a}. The characterization
of Lagrangian higher order semisprays in terms of a closed 2-form
appears also in~\cite{anderso92}.

Based on this result we can obtain the following reformulation in
terms of semi-basic 1-forms of the known Helmholtz conditions \cite[Lemma~4.2, Lemma~4.3,  Theorem~4.5, Theorem~5.1]{B-C-2010}.

\begin{theorem}
A semispray $S$ is a Lagrangian vector field if and only if there
exists a semi-basic $1$-form $\theta\in\Lambda_{v}^{1}$, with $\operatorname{rank}(d\theta)=2n$
on $J^{1}\pi$, such that \begin{gather}
d_{J}\theta=0,\qquad  d_{h}\theta=0.\label{eq:P}
\end{gather}
\end{theorem}
\begin{proof}
In order to make this paper self contained, we give a direct proof
of this theorem.

Suppose that $S$ is a Lagrangian semispray. It results that there
exists a regular Lagrangian~$L$ on~$J^{1}\pi$ with $\mathcal{L}_{S}\theta_{L}=dL$,
or equivalently $i_{S}d\theta_{L}=0$, where $\theta_{L}=Ldt+d_{J}L$
is its Poincar\'e 1-form. Evidently $\theta_{L}$ is a semi-basic 1-form
with $\operatorname{rank}(d\theta_{L})=2n$ on $J^{1}\pi$. We will
prove that $d_{J}\theta_{L}=d_{h}\theta_{L}=0$.

Indeed, $d_{J}\theta_{L}=d_{J}L\wedge dt+Ld_{J}dt+d_{J}^{2}L=i_{J}dL\wedge dt-d_{J\wedge dt}L=i_{J}dL\wedge dt-i_{J\wedge dt}dL=0$.

From the formula $i_{J}\mathcal{L}_{S}-\mathcal{L}_{S}i_{J}=i_{\Gamma-S\otimes dt}$
and $i_{J}\mathcal{L}_{S}d\theta_{L}=0$ we obtain $\mathcal{L}_{S}d_{J}\theta_{L}+i_{\Gamma}d\theta_{L}-i_{S\otimes dt}d\theta_{L}=0$.
We also compute $i_{\Gamma}d\theta_{L}=i_{2h-\operatorname{Id}}d\theta_{L}=2i_{h}d\theta_{L}-2d\theta_{L}=2d_{h}\theta_{L}$.
Therefore $2d_{h}\theta_{L}=i_{S\otimes dt}d\theta_{L}=-i_{S}d\theta_{L}\wedge dt=0$.

Conversely, suppose that there exists a semi-basic $1$-form $\theta\in\Lambda_{v}^{1}$,
with $\operatorname{rank}(d\theta)=2n$ on~$J^{1}\pi$, such that
$d_{J}\theta=0$, $d_{h}\theta=0$. In order to prove that $S$ is
a Lagrangian vector f\/ield, we will f\/irst show that~$\mathcal{L}_{S}\theta=d(i_{S}\theta)$.

The hypothesis $d_{J}\theta=0$ implies $\theta=(i_{S}\theta)dt+d_{J}(i_{S}\theta).$
Indeed, $d_{J}i_{S}+i_{S}d_{J}=\mathcal{L}_{JS}-i_{[S,J]}=i_{h-S\otimes dt-v}\Rightarrow d_{J}(i_{S}\theta)=i_{h}\theta-(i_{S}\theta)dt-i_{v}\theta=\theta-(i_{S}\theta)dt.$

Next, from $d_{h}i_{S}+i_{S}d_{h}=\mathcal{L}_{hS}-i_{[S,h]}$ and
$d_{h}\theta=0$ it results that $d_{h}i_{S}\theta=\mathcal{L}_{S}\theta-i_{\mathbb{F}+J+\Phi}\theta=\mathcal{L}_{S}\theta-i_{\mathbb{F}}\theta$.

{\sloppy From $i_{\mathbb{F}}\theta=i_{\mathbb{F}}\left((i_{S}\theta)dt+d_{J}(i_{S}\theta)\right)=i_{\mathbb{F}}\left(d_{J}(i_{S}\theta)\right)$
and $i_{\mathbb{F}}d_{J}-d_{J}i_{\mathbb{F}}=d_{J\circ\mathbb{F}}-i_{[\mathbb{F},J]}\Rightarrow $ $ i_{\mathbb{F}}\left(d_{J}(i_{S}\theta)\right)=d_{v}(i_{S}\theta)$.
It results that $\mathcal{L}_{S}\theta=d_{h}(i_{S}\theta)+d_{v}(i_{S}\theta)=d(i_{S}\theta)$.

}

Consider $L=i_{S}\theta$. Then $\theta$ is the Poincar\'e--Cartan 1-form
of $L$ and $\mathcal{L}_{S}\theta_{L}=dL$. From $\operatorname{rank}(d\theta)=2n$
on $J^{1}\pi$ it results that $L$ is a regular Lagrangian and $S$
is a Lagrangian vector f\/ield.
\end{proof}

In the next section we discuss the formal integrability of these Helmholtz
conditions using two suf\/f\/icient conditions provided by Cartan--K\"ahler
theorem.

\section{Formal integrability for the nonautonomus inverse problem\\ of the
calculus of variations}\label{section3}

In order to study the integrability conditions of the set of dif\/ferential
equations \eqref{eq:P}, we associate to it a linear  partial dif\/ferential
operator and study its formal integrability, using Spencer's technique.
The approach in this work follows the one developed in \cite{B-M_arxiv}
for studying the projective metrizability problem for autonomous sprays.
For the basic notions of formal integrability theory of linear partial
dif\/ferential operators see \cite{B-M_arxiv,grifone00}.

Consider $T_{v}^{*}$ the vector bundle of semi-basic $1$-forms on
$J^{1}\pi$ and $\Lambda_{v}^{1}$ the module of sections of~$T_{v}^{*}$.
For $\theta\in\Lambda_{v}^{1}$ and $k\geq1$ we denote by $j_{u}^{k}\theta$
the $k$th order jet of $\theta$ at the base point $u$ in $J^{1}\pi$.
The bundle of $k$th order jets of sections of $T_{v}^{*}$ is denoted
by $J^{k}T_{v}^{*}$. The projection $\pi_{0}:J^{k}T_{v}^{*}\rightarrow J^{1}\pi$
is def\/ined by $\pi_{0}(j_{u}^{k}\theta)=u$. If $l>k$, one def\/ines
the projections $\pi_{k}$ as follows: $\pi_{k}(j_{u}^{l}\theta)=j_{u}^{k}\theta$
and $J^{l}T_{v}^{*}$ is also a f\/ibred manifold over $J^{k}T_{v}^{*}$.

If $f_{1},\dots,f_{k}\in C^{\infty}(J^{1}\pi)$ are functions vanishing
at $u\in J^{1}\pi$ and $\theta\in\Lambda_{v}^{1}$, we def\/ine $\epsilon: \mathcal{S}^{k}T^{*}\otimes T_{v}^{*}\overset{}{\longrightarrow}J^{k}T_{v}^{*}$
by $\epsilon(df_{1}\odot\cdots\odot df_{k}\otimes\theta)_{u}=j_{u}^{k}(f_{1}\cdots f_{k}\theta)$,
where $\odot$ is the symmetric product. Then the sequence
\[
0\overset{}{\longrightarrow}\mathcal{S}^{k}T^{*}\otimes T_{v}^{*}\overset{\epsilon}{\longrightarrow}J^{k}T_{v}^{*}\overset{\pi_{k-1}}{\longrightarrow}J^{k-1}T_{v}^{*}\overset{}{\longrightarrow}0\]
 is exact.

Consider the \emph{linear partial differential operator} of order
one\begin{gather}
P: \ \Lambda_{v}^{1}   \rightarrow   \Lambda_{v}^{2}\oplus\Lambda_{v}^{2},\qquad
P   =  \left(d_{J},  d_{h}\right).\label{eq:PDE}
\end{gather}

Remark that $P(\theta)$ can be expressed in terms of f\/irst-order
jets of $\theta$, for any $\theta\in\Lambda_{v}^{1}$, and therefore
it induces a morphism between vector bundles:
\begin{gather*}
p^{0}(P): \  J^{1}T_{v}^{*}   \rightarrow   \Lambda^{2}T_{v}^{*}\oplus\Lambda^{2}T_{v}^{*},\qquad
p^{0}(P)(j_{u}^{1}\theta)   =    P(\theta)_{u},\quad\forall\,\theta\in\Lambda_{v}^{1}.
\end{gather*}

We also consider the \emph{$l$th order jet prolongations} of the
dif\/ferential operator $P$, $l\geq1$, which will be identif\/ied with
the morphisms of vector bundles over $M$,
\begin{gather*}
p^{l}(P): \  J^{l+1}T_{v}^{*}   \rightarrow   J^{l}\left(\Lambda^{2}T_{v}^{*}\oplus\Lambda^{2}T_{v}^{*}\right),\qquad
p^{l}(P)\big(j_{u}^{l+1}\theta\big)   =   j_{u}^{l}\left(P(\theta)\right),\quad \forall\,\theta\in\Lambda_{v}^{1}.
\end{gather*}

Remark that for a semi-basic 1-form $\theta=\theta_{\alpha}\delta x^{\alpha}$,
its f\/irst-order jet $j^{1}\theta=\frac{\delta\theta_{\alpha}}{\delta x^{\beta}}\delta x^{\beta}\otimes\delta x^{\alpha}+\frac{\partial\theta_{\alpha}}{\partial y^{i}}\delta y^{i}\otimes\delta x^{\alpha}$
determines the local coordinates $\left(x^{\alpha},y^{i},\theta_{\alpha},\theta_{\alpha\beta},\theta_{\alpha\underline{i}}\right)$
on $J^{1}T_{v}^{*}$. In this work all contravariant or covariant
indices, related to vertical components of tensor f\/ields will be underlined.

Consider $\theta=\theta_{\alpha}\delta x^{\alpha}$, a semi-basic
$1$-form on $J^{1}\pi$. Then
\begin{gather*}
d\theta   =   \left(\frac{\partial\theta_{i}}{\partial t}-\theta_{j}N_{i}^{j}-\frac{\delta\theta_{0}}{\delta x^{i}}\right)\delta x^{0}\wedge\delta x^{i}+\left(\theta_{i}-\frac{\partial\theta_{0}}{\partial y^{i}}\right)\delta x^{0}\wedge\delta y^{i}\\
\hphantom{d\theta   =}{}
+\frac{1}{2}\left(\frac{\delta\theta_{j}}{\delta x^{i}}-\frac{\delta\theta_{i}}{\delta x^{j}}\right)\delta x^{i}\wedge\delta x^{j}
   +\left(\frac{\partial\theta_{j}}{\partial y^{i}}\right)\delta y^{i}\wedge\delta x^{j},\\
d_{J}\theta   =   \left(\theta_{i}-\frac{\partial\theta_{0}}{\partial y^{i}}\right)\delta x^{0}\wedge\delta y^{i}+\frac{1}{2}\left(\frac{\partial\theta_{i}}{\partial y^{j}}-\frac{\partial\theta_{j}}{\partial y^{i}}\right)\delta x^{j}\wedge\delta x^{i},\\
d_{h}\theta   =   \left(\frac{\partial\theta_{i}}{\partial t}-\theta_{j}N_{i}^{j}-\frac{\delta\theta_{0}}{\delta x^{i}}\right)\delta x^{0}\wedge\delta x^{i}+\frac{1}{2}\left(\frac{\delta\theta_{i}}{\delta x^{j}}-\frac{\delta\theta_{j}}{\delta x^{i}}\right)\delta x^{j}\wedge\delta x^{i}.
\end{gather*}

Using these formulae we obtain
\begin{gather*}
p^{0}(P)\left(j^{1}\theta\right)=
\left(\left(\theta_{i}-\frac{\partial\theta_{0}}{\partial y^{i}}\right)\delta x^{0}\wedge\delta y^{i}+\frac{1}{2}\left(\frac{\partial\theta_{i}}{\partial y^{j}}-\frac{\partial\theta_{j}}{\partial y^{i}}\right)\delta x^{j}\wedge\delta x^{i} , \right.\\
\left.\phantom{p^{0}(P)\left(j^{1}\theta\right)= }{}
\left(\frac{\partial\theta_{i}}{\partial t}-\theta_{j}N_{i}^{j}-\frac{\delta\theta_{0}}{\delta x^{i}}\right)\delta x^{0}\wedge\delta x^{i}+\frac{1}{2}\left(\frac{\delta\theta_{i}}{\delta x^{j}}-\frac{\delta\theta_{j}}{\delta x^{i}}\right)\delta x^{j}\wedge\delta x^{i}\right).
\end{gather*}

The \emph{symbol} of $P$ is the vector bundle morphism $\sigma^{1}(P):  T^{*}\otimes T_{v}^{*}\rightarrow\Lambda^{2}T_{v}^{*}\oplus\Lambda^{2}T_{v}^{*}$
def\/ined by the f\/irst-order terms of $p^{0}(P)$. More exactly, $\sigma^{1}(P)=p^{0}(P)\circ\epsilon$.

For $A\in T^{*}\otimes T_{v}^{*}$, $A=A_{\alpha\beta}\delta x^{\alpha}\otimes\delta x^{\beta}+A_{\underline{i}\beta}\delta y^{i}\otimes\delta x^{\beta}$,
we compute
\begin{gather*}
\sigma^{1}(d_{J})A   =   -A_{\underline{i}0}\delta x^{0}\wedge\delta x^{i}+\frac{1}{2}\big(A_{\underline{j}i}-A_{\underline{i}j}\big)\delta x^{j}\wedge\delta x^{i},%\label{eq:sigma_1_dJ}
\\
\sigma^{1}(d_{h})A   =   \left(A_{0i}-A_{i0}\right)\delta x^{0}\wedge\delta x^{i}+\frac{1}{2}\left(A_{ji}-A_{ij}\right)\delta x^{j}\wedge\delta x^{i}
%\label{eq:sigma_1_dh}
\end{gather*}
 and hence \begin{gather*}
\sigma^{1}(P)A   =   \left(\tau_{J}A , \tau_{h}A\right),\qquad %\label{eq:sigma_1P}\\
\left(\tau_{J}A\right)(X,Y)   =   A(JX,Y)-A(JY,X), \\
\left(\tau_{h}A\right)(X,Y)   =   A(hX,Y)-A(hY,X),
 \end{gather*}
 for $X,Y\in\mathfrak{X}(J^{1}\pi)$. In the above formulae $\tau_{J}$, $\tau_{L}$
are \emph{alternating operators}~\cite{B-M_arxiv}.

\begin{remark}
The alternating operators are def\/ined in general as follows. For $K\in\Psi^{k}(J^{1}\pi)$,
a~vector-valued $k$-form, we consider $\tau_{K}:\Psi^{1}(J^{1}\pi)\otimes\Psi^{l}(J^{1}\pi)\to\Psi^{l+k}(J^{1}\pi)$,
\begin{gather}
  (\tau_{K}B)(X_{1},\dots ,X_{l+k})=
  \frac{1}{l!k!}\sum_{\sigma\in S_{l+k}}\varepsilon(\sigma)B(K(X_{\sigma(1)},\dots ,X_{\sigma(k)}),X_{\sigma(k+1)},\dots ,X_{\sigma(k+l)}),\label{taull}  \end{gather}
 where $X_{1},\dots ,X_{l+k}\in{\mathfrak{X}}(J^{1}\pi)$ and $S_{l+k}$
is the permutation group of $\{1,\dots ,l+k\}$. The restriction of $\tau_{K}$
to $\Psi^{l+1}(J^{1}\pi)$, is a derivation of degree $(k-1)$ and
it coincides with the \emph{inner product} $i_{K}$.
\end{remark}

The \emph{first-order prolongation of the symbol} of $P$ is the vector
bundle morphism
$\sigma^{2}(P):  S^{2}T^{*}\otimes T_{v}^{*}\rightarrow T^{*}\otimes\left(\Lambda^{2}T_{v}^{*}\oplus\Lambda^{2}T_{v}^{*}\right)$
that verif\/ies \[
i_{X}\left(\sigma^{2}(P)B\right)=\sigma^{1}(P)\left(i_{X}B\right),\qquad \forall\, B\in S^{2}T^{*}\otimes T_{v}^{*},\quad \forall\, X\in\mathfrak{X}(J^{1}\pi).
\]

Therefore \begin{gather*}
\sigma^{2}(P)B=\left(\sigma^{2}(d_{J})B , \sigma^{2}(d_{h})B\right),\\
\left(\sigma^{2}(d_{J})B\right)(X,Y,Z)   =   B(X,JY,Z)-B(X,JZ,Y), \\ %\label{eq:sigma_2_dJ}\\
\left(\sigma^{2}(d_{h})B\right)(X,Y,Z)   =   B(X,hY,Z)-B(X,hZ,Y). %\label{eq:sigma_2_dh}
\end{gather*}

In local coordinates we obtain the following formulae.

If $B\in S^{2}T^{*}\otimes T_{v}^{*}$, then it has the local decomposition
\begin{gather}
B   =   B_{\alpha\beta\gamma}\delta x^{\alpha}\otimes\delta x^{\beta}\otimes\delta x^{\gamma}+B_{\underline{i}\alpha\beta}\delta y^{i}\otimes\delta x^{\alpha}\otimes\delta x^{\beta}\nonumber \\
\hphantom{B   =}{} +B_{\alpha\underline{i}\beta}\delta x^{\alpha}\otimes\delta y^{i}\otimes\delta x^{\beta}+B_{\underline{ij}\alpha}\delta y^{i}\otimes\delta y^{j}\otimes\delta x^{\alpha},\label{eq:localB}
 \end{gather}
 with \begin{gather}
B_{\alpha\beta\gamma}=B_{\beta\alpha\gamma},\qquad B_{\underline{i}j\alpha}=B_{i\underline{j}\alpha},\qquad B_{\underline{i}0\alpha}=B_{0\underline{i}\alpha},\qquad B_{\underline{ij}\alpha}=B_{\underline{ji}\alpha}.\label{eq:condB}
\end{gather}

The f\/irst-order prolongation of the symbol of $P$ is given by
\begin{gather*}
\sigma^{2}(d_{J})B   =   B_{\alpha\underline{i}0}\delta x^{\alpha}\otimes\delta x^{i}\wedge\delta x^{0}+B_{\underline{ij}0}\delta y^{i}\otimes\delta x^{j}\wedge\delta x^{0}\nonumber \\
\hphantom{\sigma^{2}(d_{J})B   =}{}
 +\frac{1}{2}\big(B_{\alpha\underline{i}j}-B_{\alpha\underline{j}i}\big)\delta x^{\alpha}\otimes\delta x^{i}\wedge\delta x^{j}+\frac{1}{2}\big(B_{\underline{ij}k}-B_{\underline{ik}j}\big)\delta y^{i}\otimes\delta x^{j}\wedge\delta x^{k},\\ %\label{eq:localsigma2dJ}\\
\sigma^{2}(d_{h})B   =   \frac{1}{2}\left(B_{\alpha\beta\gamma}-B_{\alpha\gamma\beta}\right)\delta x^{\alpha}\otimes\delta x^{\beta}\wedge\delta x^{\gamma}+\frac{1}{2}\left(B_{\underline{i}\alpha\beta}-B_{\underline{i}\beta\alpha}\right)\delta y^{i}\otimes\delta x^{\alpha}\wedge\delta x^{\beta}. \end{gather*}

For each $u\in J^{1}\pi$, we consider
\begin{gather*}
g_{u}^{k}(P)   =   \operatorname{Ker}\sigma_{u}^{k}(P),\qquad k\in\{1,2\},\\
g_{u}^{1}(P)_{e_{1}\dots e_{j}}   =   \{A\in g_{u}^{1}(P)|i_{e_{1}}A=\cdots=i_{e_{j}}A=0\},\qquad j\in\{1,\dots ,n\},
\end{gather*}
 where $\{e_{1},\dots ,e_{n}\}$ is a basis of $T_{u}(J^{1}\pi)$. Such
a basis is called \emph{quasi-regular} if it satisf\/ies
\begin{gather*}
\dim g_{u}^{2}(P)=\dim g_{u}^{1}(P)+\sum_{j=1}^{n}\dim g_{u}^{1}(P)_{e_{1}\dots e_{j}}.%\label{quasir}
\end{gather*}

\begin{definition}
The symbol $\sigma^{1}(\mathcal{P})$ is called \emph{involutive}
at $u$ in $J^{1}\pi$ if there exists a quasi-regular basis of $T_{u}J^{1}\pi$.
\end{definition}

A f\/irst-order jet $j_{u}^{1}\theta\in J^{1}T_{v}^{*}$ is \emph{a~first-order formal solution of~$P$ at~$u$} in $J^{1}\pi$ if $p^{0}(P)(\theta)_{u}=0$.

For $l\geq1$, a $(1+l)$th order jet $j_{u}^{1+l}\theta\in J_{u}^{1+l}T_{v}^{*}$
is \emph{a $(1+l)$th order formal solution of $P$ at $u$} in $J^{1}\pi$
if $p^{l}(P)(\theta)_{u}=0$.

For any $l\geq0$, consider $R_{u}^{1+l}(P)=\ker p_{u}^{l}(P)$ the
space of \emph{$(1+l)$th order formal solutions} of~$P$ at~$u$.
We denote also $\bar{\pi}_{l,u} :  R_{u}^{1+l}(P)\rightarrow R_{u}^{l}(P)$
the restriction of $\pi_{l,u} :  J_{u}^{1+l}(T_{v}^{*})\rightarrow J_{u}^{l}(T_{v}^{*})$
to $R_{u}^{1+l}(P)$.

\begin{definition}
The partial dif\/ferential operator $P$ is called \emph{formally integrable}
at $u$ in $J^{1}\pi$ if $R^{1+l}(P)=\bigcup_{u\in J^{1}\pi}R_{u}^{1+l}(P)$
is a vector bundle over $J^{1}\pi$, for all $l\geq0$, and the map
$\bar{\pi}_{l,u}:R_{u}^{1+l}(P)\rightarrow R_{u}^{l}(P)$ is onto
for all $l\geq1$.
\end{definition}

The f\/ibred submanifold $R^{1}(P)$ of $\pi_{0}:J_{u}^{1}(T_{v}^{*})\rightarrow J^{1}\pi$
is called \emph{the partial differential equation corresponding to
the first-order PDO}~$P$. A solution of the operator $P$ on an open
set $U\subset J^{1}\pi$ is a section $\theta\in\Lambda_{v}^{1}$
def\/ined on $U$ such that $P\theta=0\Leftrightarrow p^{0}(P)(j_{u}^{1}\theta)=0$, $\forall\, u\in U$.

The Cartan--K\"ahler theorem \cite{grifone00} takes the following form
for the particular case of f\/irst-order PDO.

\begin{theorem}
Let $P$ be a first-order linear partial differential operator with
$g^{2}(P)$ a vector bundle over $R^{1}(P)$. If $\overline{\pi}_{1}:R^{2}(P)\to R^{1}(P)$
is onto and the symbol $\sigma^{1}(P)$ is involutive, then $P$ is
formally integrable.
\end{theorem}

\subsection[The involutivity of the symbol of $P$]{The involutivity of the symbol of $\boldsymbol{P}$}

In this subsection we prove that the operator $P$ satisf\/ies one of
the two suf\/f\/icient conditions for formal integrability, provided by
Cartan--K\"ahler theorem: the involutivity of the symbol~$\sigma^{1}(P)$.

\begin{theorem}\label{thm:-involutivity}
The symbol $\sigma^{1}(P)$ of the PDO
$P=(d_{J},d_{h})$ is involutive.
\end{theorem}

\begin{proof}
First we determine $g^{1}(P)=\left\{   A\in T^{*}\otimes T_{v}^{*}\,\vert\,\sigma^{1}(P)A=0\right\} $,
and compute the dimension of its f\/ibers. We obtain
\begin{gather*}
g_{u}^{1}(P)   =   \big\{  A=A_{\alpha\beta}\delta x^{\alpha}\otimes\delta x^{\beta}+A_{\underline{i}\beta}\delta y^{i}\otimes\delta x^{\beta}\,\,\vert\,\, A_{\underline{i}0}=0,\, A_{\underline{j}i}=A_{\underline{i}j},\, A_{\alpha\beta}=A_{\beta\alpha}\big\} .
\end{gather*}
 From $A_{\underline{i}0}=0$ and $A_{\underline{j}i}=A_{\underline{i}j}$
it results that $A_{\underline{i}j}$ contribute with $n(n+1)/2$
components to the dimension of $g_{u}^{1}(P)$, and from $A_{\alpha\beta}=A_{\beta\alpha}$
it follows that $A_{\alpha\beta}$ contribute with $(n+1)(n+2)/2$
components to the dimension of $g_{u}^{1}(P)$. So
\[
\dim g_{u}^{1}(P)=\frac{n(n+1)}{2}+\frac{(n+1)(n+2)}{2}=(n+1)^{2}.
\]

Next we determine $g^{2}(P)=\left\{   B\in S^{2}T^{*}\otimes T_{v}^{*}\,\vert\,\sigma^{2}(P)B=0\right\} $.

If $B\in S^{2}T^{*}\otimes T_{v}^{*}$ has the local components (\ref{eq:localB}),
then $B\in g^{2}(P)$ if and only if the following relations are satisf\/ied:
\begin{gather}
  B_{\alpha\underline{i}0}=0,  \qquad    B_{\underline{ij}0}=0, \qquad
  B_{\alpha\underline{i}j}=B_{\alpha\underline{j}i}, \nonumber \\
  B_{\underline{ij}k}=B_{\underline{ik}j},\qquad
  B_{\alpha\beta\gamma}=B_{\alpha\gamma\beta}, \qquad B_{\underline{i}\alpha\beta}=B_{\underline{i}\beta\alpha}.\label{eq:eq:cond_B_in_g2}
  \end{gather}

From the relations \eqref{eq:condB} and \eqref{eq:eq:cond_B_in_g2}
it results that $B\in g^{2}(P)$ if and only if its local components
$  B_{\alpha\beta\gamma}$, $B_{\underline{i}jk}$, $B_{i\underline{j}k}$, $B_{\underline{ij}k}$
are totally symmetric and the rest are vanishing. Therefore $B_{\alpha\beta\gamma}$
contribute with $(n+1)(n+2)(n+3)/6$ components to the dimension of
$g_{u}^{2}(P)$, and $B_{\underline{i}jk}$,  $B_{\underline{ij}k}$
with $n(n+1)(n+2)/6$ components each of them. It results
\[
\dim g_{u}^{2}(P)=\frac{(n+1)(n+2)(n+3)}{6}+2\frac{n(n+1)(n+2)}{6}=\frac{(n+1)^{2}(n+2)}{2}.
\]

Consider \[
\mathcal{B}=\left\{ h_{0}=S,\, h_{1},\,\dots,h_{n},\, v_{1}=Jh_{1},\,\dots,\, v_{n}=Jh_{n}\right\} \]
 a basis in $T_{u}J^{1}\pi$ with $h_{0}=S$, $h_{1},\dots,h_{n}$
horizontal vector f\/ields. For any $A\in g^{1}(P)$, we denote \[
A(h_{\alpha},h_{\beta})=a_{\alpha\beta},\qquad A(v_{i},h_{\alpha})=b_{\underline{i}\alpha}.\]

Because $A$ is semi-basic in the second argument it follows that
these are the only components of $A$.

Since $A\in\operatorname{Ker}\sigma^{1}(d_{J})$ it follows that $A_{\underline{i}0}=0$, $A_{\underline{i}j}=A_{\underline{j}i}$ $\Rightarrow$ $b_{\underline{i}0}=0$, $b_{\underline{i}j}=b_{\underline{j}i}$.
Because $A\in\operatorname{Ker}\sigma^{1}(d_{h})$ it results that
$A_{0i}=A_{i0}$, $A_{ij}=A_{ji}$ and hence $a_{0i}=a_{i0}$, $a_{ij}=a_{ji}$.

Consider $j\in \{1,\dots, n\} $ arbitrarily f\/ixed and
\begin{gather*}
\mathcal{\tilde{B}}=\big\{ e_{0}=S+h_{j}+v_{n},\, e_{1}=h_{1},\, e_{2}=h_{2}+v_{1}, \, \dots,\\
 \hphantom{\mathcal{\tilde{B}}=\big\{}{} e_{i}=h_{i}+v_{i-1},\, \dots,\, e_{n}=h_{n}+v_{n-1},\, v_{1},\, \dots,\, v_{n}\big\}
\end{gather*}
 a new basis in $T_{u}J^{1}\pi$. If we denote \[
A(e_{\alpha},e_{\beta})=\tilde{a}_{\alpha\beta},\qquad A(v_{i},e_{\alpha})=\tilde{b}_{\underline{i}\alpha},\]
 a simple computation and the fact that $A$ is semi-basic in the
second argument determine \begin{gather*}
\tilde{a}_{00}   =   a_{00}+2a_{0j}+a_{jj}+b_{\underline{n}j},\\
\tilde{a}_{ik}   =   a_{ik}+b_{\underline{i-1},k}\neq\tilde{a}_{ki}=a_{ki}+b_{\underline{k-1},i},\\
\tilde{a}_{i0}   =   a_{i0}+a_{ij}+b_{\underline{i-1},j}\neq\tilde{a}_{0i}=a_{0i}+a_{ji}+b_{\underline{n}i},\\
\tilde{b}_{\underline{i}k}   =b_{\underline{i}k}   =\tilde{b}_{\underline{k}i},\qquad
\tilde{b}_{\underline{i}0}   =   b_{\underline{i}j}.
\end{gather*}

It can be seen that all the independent components of $A$ in the
basis $\mathcal{B}$ can be obtained from the components of $A$ in
the basis $\mathcal{\tilde{B}}$, and hence we can use the later for
determining the dimensions of $(g_{u}^{1})_{e_{0} e_{1}\dots e_{k}}$.

If $A\in(g_{u}^{1})_{e_{0}}$ it results $\tilde{a}_{0\alpha}=0$,
so using this new basis we impose $n+1$ supplementary independent
restrictions. It follows that $\dim(g_{u}^{1})_{e_{0}}=(n+1)^{2}-(n+1)=(n+1)n$.

If $A\in(g_{u}^{1})_{e_{0} e_{1}}$ it results that together with
the previous restrictions we impose also $\tilde{a}_{1\alpha}=0$,
so another independent $n+1$ restrictions. Hence $\dim (g_{u}^{1})_{e_{0} e_{1}}=(n+1)n-(n+1)=(n+1)(n-1)$.

In general $\dim (g_{u}^{1})_{e_{0} e_{1} \dots e_{k}}=(n+1)(n-k)$, $\forall\, k\in \{1,\dots,n\}$.
Hence $\dim (g_{u}^{1})_{e_{0} e_{1} \dots e_{n}}=0\Rightarrow\dim (g_{u}^{1})_{e_{0} \dots e_{n} v_{1} \dots v_{k}}=0$, $\forall\, k\in \{1,\dots,n\}$,
\begin{gather*}
\dim (g_{u}^{1})+\sum_{k=0}^{n}\dim (g_{u}^{1})_{e_{0} e_{1} \dots e_{k}}
   =   (n+1)^{2}+(n+1)[n+n-1+\cdots+1]\\
\hphantom{\dim (g_{u}^{1})+\sum_{k=0}^{n}\dim (g_{u}^{1})_{e_{0} e_{1} \dots e_{k}}}{}   =   \frac{(n+1)^{2}(n+2)}{2}=\dim g_{u}^{2}(P).
   \end{gather*}

We proved that $\tilde{B}$ is a a quasi-regular basis, hence the
symbol $\sigma^{1}(P)$ is involutive.
\end{proof}

\subsection{First obstruction to the inverse problem}

In this subsection we determine necessary and suf\/f\/icient conditions
for $\bar{\pi}_{1}$ to be onto. We will obtain only one obstruction
for the integrability of the operator $P$. The obstruction is due
to the curvature tensor of the nonlinear connection induced by the
semispray.

\begin{theorem}
\label{thm:A-first-order}A first-order formal solution $\theta\in\Lambda_{v}^{1}$
of the system $d_{J}\theta=0$, $d_{h}\theta=0$ can be lifted into
a second-order solution, which means that $\bar{\pi}_{1} :  R^{2}(P)\rightarrow R^{1}(P)$
is onto, if and only if \[
d_{R}\theta=0,\]
 where $R$ is the curvature tensor \eqref{eq:19}.
\end{theorem}
\begin{proof}
We use a known result from \cite[Proposition~1.1]{grifone00}.

If $\mathcal{K}$ is the cokernel of $\sigma^{2}(P)$,
\begin{gather*}
\mathcal{K}=\frac{T^{*}\otimes\left(\Lambda^{2}T_{v}^{*}\oplus\Lambda^{2}T_{v}^{*}\right)}{{\normalcolor \operatorname{Im}}\sigma^{2}(P)},%\label{eq:K}
\end{gather*}
 there exists a morphism $\varphi :R^{1}(P)\rightarrow\mathcal{K}$
such that the sequence \[
R^{2}(P)\overset{\bar{\pi}_{1}}{\longrightarrow}R^{1}(P)\overset{\varphi}{\longrightarrow}\mathcal{K}\]
 is exact. In particular $\bar{\pi}_{1}$ is onto if and only if $\varphi=0$.

After def\/ining $\varphi$, we will prove that for $\theta\in\Lambda_{v}^{1}$,
with $j_{u}^{1}\theta\in R_{u}^{1}(P)$, a f\/irst-order formal solution
of $P$ at $u\in J^{1}\pi$, we have that $\varphi_{u}\theta=0$ if
and only if $\left(d_{R}\theta\right)_{u}=0$.

The construction of the morphism $\varphi$ is represented in the
next diagram by dashed arrows. We denote $F=\Lambda^{2}T_{v}^{*}\oplus\Lambda^{2}T_{v}^{*}$.
%$$  {\diagram & 0 \dto & 0 \dto & 0 \dto \\ 0\rto & g^2(P) \rto \dto & S^{2}T^*\otimes T^*_v \rto^{\sigma^2(P)} \dto^{\varepsilon} & T^* \otimes F %\rto \rdashed<1ex>|>{\tip}^{\tau} \dto^{\varepsilon} &  \it{K} \rto & 0 \\ 0 \rto & R^{2}(P) \rto^{i} \dto^{\overline {\pi}_{1}} & J^{2} T^*_v \rto %\dto^{\pi_{1}} \rdashed<1ex>|>{\tip}^{p^1(P)} & J^1 F \dto^{\pi_0} \udashed<1ex>|>{\tip}^{p^0(\nabla)} \\ 0 \rto & R^1(P) \rto\rdashed<1ex>|>{\tip}^i %& J^1T^*_v \rto^{p^0(P)} \udashed<1ex>|>{\tip} \dto & F \dto \\ & & 0 & 0 \enddiagram}
%$$
$$
\includegraphics{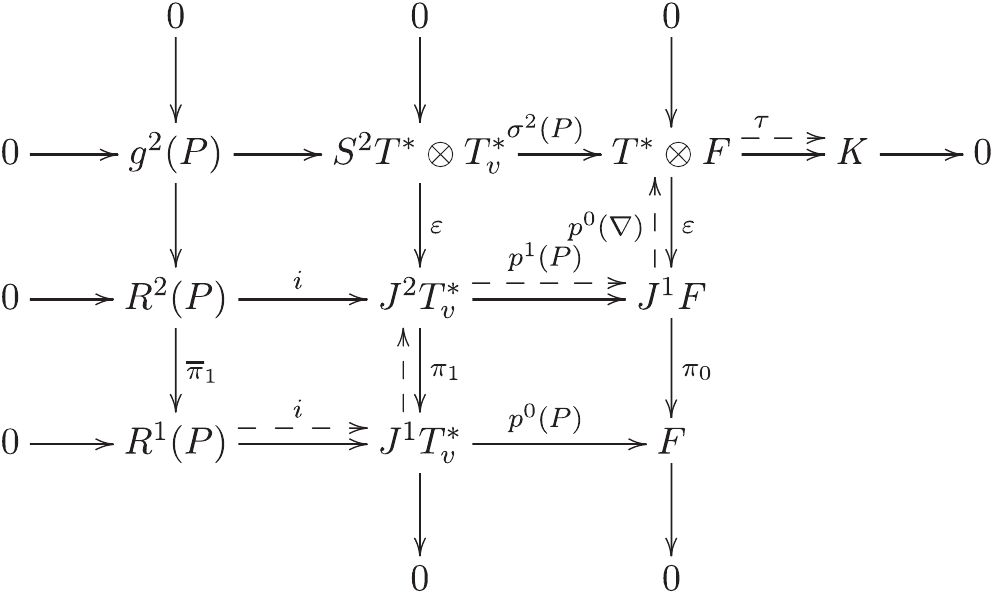}
$$

Remark that $\dim T^{*}=2n+1$, $\dim T_{v}^{*}=n+1$, $\dim\Lambda^{2}T_{v}^{*}=\frac{(n+1)n}{2}$, $\dim\mathcal{S}^{2}T^{*}=\frac{(2n+1)(2n+2)}{2}$,
$\dim T^{*}\otimes\left(\Lambda^{2}T_{v}^{*}\oplus\Lambda^{2}T_{v}^{*}\right)=(2n+1)n(n+1)$.
Therefore
\begin{gather*}
\dim\mathcal{K}   =   \dim\left[T^{*}\otimes\left(\Lambda^{2}T_{v}^{*}\oplus\Lambda^{2}T_{v}^{*}\right)\right]-\dim\left(\operatorname{Im}\sigma^{2}(P)\right)\\
\hphantom{\dim\mathcal{K}}{} =   \dim\left[T^{*}\otimes\left(\Lambda^{2}T_{v}^{*}\oplus\Lambda^{2}T_{v}^{*}\right)\right]-\left[\dim\left(\mathcal{S}^{2}T^{*}\otimes T_{v}^{*}\right)-\dim\left(\ker\sigma^{2}(P)\right)\right]\\
\hphantom{\dim\mathcal{K}}{}
 =  \frac{(n-1)n(n+1)}{2}=3 \begin{pmatrix}
n+1\\
3\end{pmatrix}.\end{gather*}

It results from this that \[
\mathcal{K}\simeq\oplus^{(3)}\Lambda^{3}T_{v}^{*}.\]

Next we def\/ine $\tau:  T^{*}\otimes\left(\Lambda^{2}T_{v}^{*}\oplus\Lambda^{2}T_{v}^{*}\right)\rightarrow\oplus^{(3)}\Lambda^{3}T_{v}^{*}$
such as the next sequence is exact:
 \begin{gather}
0\rightarrow g^{2}(P)\overset{i}{\longrightarrow}\mathcal{S}^{2}T^{*}\otimes T_{v}^{*}\overset{\sigma^{2}(P)}{\longrightarrow}T^{*}
\otimes\left(\Lambda^{2}T_{v}^{*}\oplus\Lambda^{2}T_{v}^{*}\right)\overset{\tau}{\longrightarrow}\oplus^{(3)}\Lambda^{3}T_{v}^{*}\rightarrow0.
\label{eq:exact_tau}\end{gather}

For $B_{1},  B_{2}\in T^{*}\otimes\Lambda^{2}T_{v}^{*}$, we def\/ine
$\tau(B_{1},B_{2})=\left(\tau_{1}(B_{1},B_{2}) , \tau_{2}(B_{1},B_{2}) , \tau_{3}(B_{1},B_{2})\right)$,
where $\tau_{i}:  T^{*}\otimes\left(\Lambda^{2}T_{v}^{*}\oplus\Lambda^{2}T_{v}^{*}\right)\rightarrow\Lambda^{3}T_{v}^{*}$, $ i\in \{1,2,3\}$,
are given by
\begin{gather*}
\tau_{1}(B_{1},B_{2})=\tau_{J}B_{1},\qquad
\tau_{2}(B_{1},B_{2})=\tau_{h}B_{2},\qquad \tau_{3}(B_{1},B_{2})=\tau_{h}B_{1}+\tau_{J}B_{2}.%\label{eq:tau_123}
\end{gather*}

Using the def\/inition (\ref{taull}) of the alternating operators $\tau_{J}$, $\tau_{h}$,
we prove that $\tau\circ\sigma^{2}(P)=0$.

Indeed, using that any $B\in\mathcal{S}^{2}T^{*}\otimes T_{v}^{*}$
is symmetric in the f\/irst two arguments, it follows that $\left(\tau\circ\sigma^{2}(P)\right)(B)=\tau\left(\sigma^{2}(d_{J})B, \sigma^{2}(d_{h})B\right)=\left(\tau_{J}\sigma^{2}(d_{J})B, \tau_{h}\sigma^{2}(d_{h})B, \tau_{h}\sigma^{2}(d_{J})B+\tau_{J}\sigma^{2}(d_{h})B\right)=0$, $\forall\, B\in\mathcal{S}^{2}T^{*}\otimes T_{v}^{*}$.
For example,\begin{gather*}
\tau_{J}\left(\sigma^{2}(d_{J})B\right)(X,Y,Z)  \\
 \qquad{} =  \left[\sigma^{2}(d_{J})B\left(JX,Y,Z\right)-\sigma^{2}(d_{J})B\left(JY,X,Z\right)+\sigma^{2}(d_{J})B\left(JZ,X,Y\right)\right]\\
  \qquad{}{} =   [B(JX,JY,Z)-B(JX,JZ,Y)-B(JY,JX,Z)+B(JY,JZ,X)\\
   \qquad\quad{}{}   +B(JZ,JX,Y)-B(JZ,JY,X)]=0,\qquad \forall\, X,Y,Z\in\mathfrak{X}(J^{1}\pi).
 \end{gather*}

The relation $\tau\circ\sigma^{2}(P)=0$ implies that $\operatorname{Im}(\sigma^{2}(P))\subseteq\operatorname{Ker}\tau$.
Using that $\tau$ is onto ($\tau_{J}$,~$\tau_{h}$~are both onto)
it results that $\dim\left[\operatorname{Im}(\sigma^{2}(P))\right]=\dim\left(\operatorname{Ker}\tau\right)$
and hence $\operatorname{Im}(\sigma^{2}(P))=\operatorname{Ker}\tau$
and the sequence~(\ref{eq:exact_tau}) is exact.

The last step before def\/ining $\varphi :R^{1}(P)\rightarrow\mathcal{K}$
is to consider a linear connection $\nabla$ on $J^{1}\pi$ such that
$\nabla J=0$. It means that $\nabla$ preserve semi-basic forms and
$\nabla$ can be considered as a~connection in the f\/iber bundle $\Lambda^{2}T_{v}^{*}\oplus\Lambda^{2}T_{v}^{*}\rightarrow J^{1}\pi$.
As a f\/irst-order PDO we can identify $\nabla$ with the bundle morphism
$p^{0}(\nabla):  J^{1}\left(\Lambda^{2}T_{v}^{*}\oplus\Lambda^{2}T_{v}^{*}\right)\rightarrow T^{*}\otimes\left(\Lambda^{2}T_{v}^{*}\oplus\Lambda^{2}T_{v}^{*}\right)$.

We will also use two derivations of degree~1 introduced in~\cite{B-M_arxiv},
def\/ined by $\mathcal{D}_{J}=\tau_{J}\nabla$, $\mathcal{D}_{h}=\tau_{h}\nabla.$
Both derivations $\mathcal{D}_{J}$, $\mathcal{D}_{h}$ preserve semi-basic
forms and $d_{J}-\mathcal{D}_{J}$, $ d_{h}-\mathcal{D}_{h}$ are algebraic
derivations. It means that if $\omega\in\Lambda^{k}\left(J^{1}\pi\right)$
vanishes at some point $u\in J^{1}\pi$, then $\left(\mathcal{D}_{J}\omega\right)_{u}=\left(d_{J}\omega\right)_{u}$
and $\left(\mathcal{D}_{h}\omega\right)_{u}=\left(d_{h}\omega\right)_{u}$
\cite[Lemma~2.1]{B-M_arxiv}.

Now we are able to def\/ine $\varphi :R^{1}(P)\rightarrow\mathcal{K}$
such that the sequence \[
R^{2}(P)\overset{\bar{\pi}_{1}}{\longrightarrow}R^{1}(P)\overset{\varphi}{\longrightarrow}\mathcal{K}\]
 is exact.

Let $\theta\in\Lambda_{v}^{1}$ such that $j_{u}^{1}\theta\in R_{u}^{1}(P)\subset J_{u}^{1}T_{v}^{*}$
is a f\/irst-order formal solution of $P$ at $u\in J^{1}\pi$, which
means that $\left(d_{J}\theta\right)_{u}=\left(d_{h}\theta\right)_{u}=0$.

Consider \[
\varphi_{u}\theta=\tau_{u}\nabla P\theta=\tau_{u}\left(\nabla d_{J}\theta, \nabla d_{h}\theta\right).\]

Using the fact that $d_{J}-\tau_{J}\nabla$ and $d_{h}-\tau_{h}\nabla$
are algebraic derivations and $\left(d_{J}\theta\right)_{u}=\left(d_{h}\theta\right)_{u}=0$
it results $d_{J}\left(d_{J}\theta\right)_{u}=\tau_{J}\nabla\left(d_{J}\theta\right)_{u}$
and $d_{h}\left(d_{h}\theta\right)_{u}=\tau_{h}\nabla\left(d_{h}\theta\right)_{u}$.

We will compute the three components of the map $\varphi$. It follows
that \begin{gather*}
\tau_{1}\left(\nabla P\theta\right)_{u}   =   \tau_{1}\left(\nabla d_{J}\theta, \nabla d_{h}\theta\right)_{u}=\tau_{J}\left(\nabla d_{J}\theta\right)_{u}  \\
\hphantom{\tau_{1}\left(\nabla P\theta\right)_{u}}{}
 =\left(d_{J}^{2}\theta\right)_{u}=\frac{1}{2}\left(d_{[J,J]}\theta\right)_{u}=-\left(d_{J\wedge dt}\theta\right)_{u}=-\left(d_{J}\theta\right)_{u}\wedge dt=0, %\label{eq:tau_1nablaPtheta}
 \\
\tau_{2}\left(\nabla P\theta\right)_{u}   =   \tau_{2}\left(\nabla d_{J}\theta, \nabla d_{h}\theta\right)_{u}=\tau_{h}\left(\nabla d_{h}\theta\right)_{u}
   =   \left(d_{h}^{2}\theta\right)_{u}=\frac{1}{2}\left(d_{[h,h]}\theta\right)_{u}=\left(d_{R}\theta\right)_{u},%\label{eq:obstructiond_R}
   \end{gather*}
 where $R$ is given by (\ref{eq:19}),
\begin{gather*}
\tau_{3}\left(\nabla P\theta\right)_{u}   =   \tau_{3}\left(\nabla d_{J}\theta, \nabla d_{h}\theta\right)_{u}=\tau_{h}\left(\nabla d_{J}\theta\right)_{u}+\tau_{J}\left(\nabla d_{h}\theta\right)_{u}\nonumber \\
\hphantom{\tau_{3}\left(\nabla P\theta\right)_{u}}{}
   =   \left(d_{h}d_{J}\theta\right)_{u}+\left(d_{J}d_{h}\theta\right)_{u}=\left(d_{[h,J]}\theta\right)_{u}=0.%\label{eq:tau_3nablaPtheta}
 \end{gather*}
Hence $\varphi=0$ if and only if $d_{R}\theta=0$.
\end{proof}

\begin{remark}
Locally, $d_{R}\theta$ has the following form:
\begin{gather*}
R=\underbrace{\frac{1}{2}R_{ij}^{k}\frac{\partial}{\partial y^{k}}\otimes\delta x^{i}\wedge\delta x^{j}}_{\widetilde{R}=\frac{1}{3}[J,\Phi]}-\Phi\wedge dt \ \Rightarrow \ d_{R}\theta=d_{\widetilde{R}}\theta-d_{\Phi}\theta\wedge dt,
\\
d_{\Phi}\theta   =   R_{i}^{j}\left(\theta_{j}-\frac{\partial\theta_{0}}{\partial y^{j}}\right)dt\wedge\delta x^{i}+\frac{1}{2!}\left(\frac{\partial\theta_{j}}{\partial y^{k}}R_{i}^{k}-\frac{\partial\theta_{i}}{\partial y^{k}}R_{j}^{k}\right)\delta x^{j}\wedge\delta x^{i} \quad \Rightarrow\\
d_{R}\theta   =   \frac{1}{3!}\left(a_{il}R_{jk}^{l}+a_{jl}R_{ki}^{l}+a_{kl}R_{ij}^{l}\right)\delta x^{i}\wedge\delta x^{j}\wedge\delta x^{k}+\frac{1}{2!}\left(a_{jk}R_{i}^{k}-a_{ik}R_{j}^{k}\right)dt\wedge\delta x^{i}\wedge\delta x^{j},
\end{gather*}
where we denoted $a_{ij}=\frac{\partial\theta_{i}}{\partial y^{j}}$.

Hence $d_{R}\theta=0$ if and only if $a_{il}R_{jk}^{l}+a_{jl}R_{ki}^{l}+a_{kl}R_{ij}^{l}=0$
and $a_{jk}R_{i}^{k}-a_{ik}R_{j}^{k}=0$. The f\/irst identity represents
the algebraic Bianchi identity for the curvatures of the nonlinear
connection. The second identity is one of the classical Helmholtz
condition for the multiplier matrix~$a_{ij}$. These obstructions
appear also in~\cite{anderso92}.

It can be seen that for $n=2$ the formula of $d_{R}\theta$ becomes
\[
d_{R}\theta=\frac{1}{2!}\big(a_{jk}R_{i}^{k}-a_{ik}R_{j}^{k}\big)dt\wedge\delta x^{i}\wedge\delta x^{j}=-d_{\Phi}\theta\wedge dt.\]

Therefore, for $n=2$, the obstruction is
equivalent with $d_{\Phi}\theta\wedge dt=0$.
\end{remark}

\subsection{Classes of Lagrangian time-dependent SODE}\label{section3.3}

We present now some classes of semisprays for which the obstruction
in Theorem~\ref{thm:A-first-order} is auto\-ma\-ti\-cally satisf\/ied. Therefore
the PDO $P$ is formally integrable, and hence these semisprays will
be Lagrangians SODEs. These classes of semisprays are:
\begin{itemize}\itemsep=0pt
\item f\/lat semisprays, $R=0\Leftrightarrow\Phi=0$;
\item arbitrary semisprays on 2-dimensional manifolds;
\item isotropic semisprays, $\Phi=\lambda J$, for $\lambda$ a smooth function
on $J^{1}\pi$.
\end{itemize}
All these classes of semisprays were already studied in the articles
cited in the introduction.

In the f\/lat case, the obstruction is automatically satisf\/ied.

If $\dim M=1$ then for a semi-basic $1$-form $\theta$ on $J^{1}\pi$,
$d_{R}\theta$ is a semi-basic $3$-form on $J^{1}\pi$. Because $\dim\Lambda^{3}\left(T_{v}^{*}\right)=(n+1)n(n-1)/6$
and it is zero if $n=1$, $d_{R}\theta$ will necessarily vanish.

We consider now the last case, of isotropic semisprays.

\begin{proposition}
Any isotropic semispray is a Lagrangian second-order vector field.
\end{proposition}

\begin{proof}
Assume now that $S$ is an isotropic SODE and $\theta$ a semi-basic
1-form on $J^{1}\pi$ such that $\left(d_{J}\theta\right)_{u}=\left(d_{h}\theta\right)_{u}=0$,
for some $u\in J^{1}\pi$, \begin{gather*}
\left(d_{R}\theta\right)_{u}   =   d_{\alpha\wedge J}\theta=\alpha\wedge d_{J}\theta+(-1)^{2}d\alpha\wedge i_{J}\theta=0.
\end{gather*}

We used that $i_{J}\theta=0$, $\left(d_{J}\theta\right)_{u}=0$ and
the formula~\cite{grifone00}
\[
d_{\omega\wedge K}\pi=\omega\wedge d_{K}\pi+(-1)^{q+k}d\omega\wedge i_{K}\pi,
\]
 for $\omega$ a $q$-form on $J^{1}\pi$ and $K$ a vector-valued
$k$-form on $J^{1}\pi$.

Since $\left(d_{R}\theta\right)_{u}$ vanishes, $S$ is a Lagrangian
semispray.
\end{proof}

Next we give some simple examples of Lagrangian semisprays, corresponding
to the above general classes.

We start with the semispray expressed by the SODE
\[
\frac{d^{2}x^{1}}{dt^{2}}+f\left(t,\frac{dx^{2}}{dt}\right)   =0,\qquad
\frac{d^{2}x^{2}}{dt^{2}}+g(t)   =0, \]
with $f$ an arbitrary smooth function depending only on $t$ and
$y^{2}=\frac{dx^{2}}{dt}$, and~$g$ an arbitrary smooth function
depending only on $t$. The only possible nonvanishing local component
of the Jacobi endomorphism is
\[
R_{2}^{1}=-S\big(N_{2}^{1}\big)=-\frac{1}{2}\frac{\partial^{2}f}{\partial t\partial y^{2}}+g(t)\frac{\partial^{2}f}{\partial(y^{2})^{2}}.
\]
Hence, if $\frac{\partial^{2}f}{\partial t\partial y^{2}}=2g(t)\frac{\partial^{2}f}{\partial(y^{2})^{2}}$
the semispray is f\/lat ($\Phi=0$). This example is a generalization
of the one given by Douglas \cite[(8.14)]{douglas41}.

If $\frac{\partial^{2}f}{\partial t\partial y^{2}}\neq2g(t)\frac{\partial^{2}f}{\partial(y^{2})^{2}}$,
then the semispray is isotropic  \[
\Phi= \begin{pmatrix}
0 & -\dfrac{1}{2}\dfrac{\partial^{2}f}{\partial t\partial y^{2}}+g(t)\dfrac{\partial^{2}f}{\partial(y^{2})^{2}}\vspace{1mm}\\
0 & 0\end{pmatrix} =\left(-\frac{1}{2}\frac{\partial^{2}f}{\partial t\partial y^{2}}+g(t)\frac{\partial^{2}f}{\partial(y^{2})^{2}}\right)J.\]

Another example of isotropic (or f\/lat) semispray is the one given
by the SODE
\[
\frac{d^{2}x^{1}}{dt^{2}}+f\left(t,x^{2}\right)   =0,\qquad
\frac{d^{2}x^{2}}{dt^{2}}+g(t)   =0, \]
 with $f$ an arbitrary smooth function depending only on~$t$ and~$x^{2}$,  and~$g$ an arbitrary smooth function depending only on~$t$. All the local coef\/f\/icients of the associated nonlinear connection
are vanishing. Evidently \[
\Phi= \begin{pmatrix}
0 & \dfrac{\partial f}{\partial x^{2}}\vspace{1mm}\\
0 & 0\end{pmatrix} =\left(\frac{\partial f}{\partial x^{2}}\right)J.\]

This example was treated in \cite[(6.1)]{anderso92} and \cite[(15.4)]{douglas41}
for $f(x^{2})=-x^{2}$ and $g=0$. The f\/irst paper also presents all
the Lagrangians corresponding to the given SODE.

Consider also the semispray given by the SODE \[
\frac{d^{2}x^{1}}{dt^{2}}+2\frac{dx^{2}}{dt}   =0,\qquad
\frac{d^{2}x^{2}}{dt^{2}}-\left(\frac{dx^{2}}{dt}\right)^{2}   =0. \]

Evidently \[
\Phi= \begin{pmatrix}
0 & y^{2}\\
0 & 0\end{pmatrix},\]
hence the semispray is isotropic.

\subsection*{Acknowledgements}
The author express his thanks to Ioan Bucataru for the many interesting
discussions about the paper.

\pdfbookmark[1]{References}{ref}
\LastPageEnding

\end{document}